\documentclass{amsproc}

\usepackage{}

\newtheorem{theorem}{Theorem}[section]

\newtheorem{conjecture}[theorem]{Conjecture}

\theoremstyle{definition}

\newcommand{\Q}{\mathbb{Q}}
\newcommand{\C}{\mathbb{C}}
\newcommand{\G}{\mathbb{G}}

\newcommand{\A}{\mathbb{A}}
\newcommand{\PP}{\mathbb{P}}
\def \e{\epsilon}

\theoremstyle{remark}

\numberwithin{equation}{section}

\begin{document}

\title[$p$-adic distance of special points]{The $p$-adic distance of special points to subvarieties}


\author{Jos\'e Felipe Voloch}
\address{School of Mathematics and Statistics, University of Canterbury, Private Bag 4800, Christchurch 8140, New Zealand}
\curraddr{}
\email{felipe.voloch@canterbury.ac.nz}
\urladdr{http://www.math.canterbury.ac.nz/\~{}f.voloch}
\thanks{The author was supported by the Marsden Fund administered by the Royal
Society of New Zealand. He also thanks the Clay Foundation for partially supporting his attendance to the Tate Centennial Conference.}

\dedicatory{To the memory of John Tate}



\date{}

\begin{abstract}
In a paper of Tate and the author, we conjectured a uniform bound for the $p$-adic distance of torsion points on a semiabelian variety, not lying in a subvariety, to that subvariety. We survey the progress made on that conjecture and on similar statements in analogous situations.
\end{abstract}

\maketitle

\section{Introduction}

In our paper \cite{zbMATH00972701} with Tate we conjectured a uniform bound for the $p$-adic distance of torsion points on a semiabelian variety, not lying in a subvariety, to that subvariety. See Conjecture \ref{main} below, for a precise statement. We proved it for tori by reducing it to a statement about roots of unity. Since then, there has been a lot of progress on the conjecture and a number of variants in analogous situations. We will survey some of these results with some motivation coming from our roots of unity paradigm. 

Special points are, loosely speaking, points of particular arithmetic significance in algebraic varieties. We do not give a general definition but natural examples, which we will focus on, are torsion points in semiabelian varieties, preperiodic points in algebraic dynamical system and CM points on moduli spaces of abelian varieties, in general.

There is also a corresponding notion of special subvarieties and the Manin-Mumford conjecture and its generalizations state that the Zariski closure of a set of special points is a special subvariety. Even though, strictly speaking, Manin-Mumford type questions and $p$-adic approximation questions are unrelated, the results for both questions seem to go hand in hand.

\section{Preliminaries}
\label{sec:prelim}

Let $K$ be a field, complete with respect to a non-archimedian absolute value $|\cdot|$, with ring of integers $R$. In the cases of interest, $K$ will not be locally compact. We will be looking at certain subsets of algebraic varieties $X/K$ with the goal of showing that elements of this subset do not accummulate on any subvariety. In particular, considering zero-dimensional subvarieties, these sets will be discrete. As $X(K)$ is not locally compact, discreteness is not sufficient for our goal. First, we need to define a function that measures distance to a subvariety.

If $X/K$ is a variety and $Y$ is a subvariety of $X$, we can define a distance to $Y$ from points in $X(K)$ as follows. We choose models $\mathcal{X},\mathcal{Y}$ for $X$ and $Y$ over $R$. 
For any $\e, 0 < \e \le 1$, define $M_{\e} = \{x \in K \mid |x| < \e \}$.
So, $M_1 = M$ is the maximal ideal of $R$. Let $R_{\e} = R/M_{\e}$
for $0 < \e \le 1$. Denote by $\mathcal{X}_{\e}$
the base change of $\mathcal{X}$ to $R_{\e}$ and likewise for $\mathcal{Y},P$.
Let $d(P,Y) = \inf \{ \e > 0 \mid P_{\e} \in \mathcal{Y}_{\e}\}$, provided this set
is non-empty, otherwise set $d(P,Y) = 1$. Details of this construction are in \cite{zbMATH01130578}. If $Y$ is given locally as the zero set of functions $f_j$, then $d(P,Y)=c(P)\max\{|f_j(P)|\}$, where $c(P)$ is bounded above and below by positive constants.

\section{Roots of unity}
\label{sec:roots}

\begin{theorem}
\label{thm:TV}
\cite{zbMATH00972701}
If $a_1,\ldots,a_n \in K$, there exists $c>0$ such that, for any roots of unity $\zeta_1,\ldots,\zeta_n \in K$, either
$\sum a_i \zeta_i =0$ or $|\sum a_i \zeta_i | \ge c$.
\end{theorem}

The theorem provides a lower bound to the absolute value of a linear form in roots of unity whenever the linear form is nonzero. This can be easily upgraded to a lower bound on the absolute value of a of the evaluation of a (Laurent) polynomial in roots of unity, simply because a monomial in roots of unity is itself a root of unity.

It is natural to consider than whether it is possible to obtain a similar result when evaluating other functions, such as power series, evaluated on roots of unity. This question has been addressed by A. Neira \cite{Nei02},
who considered power series in $R[[x_1,\ldots,x_n]]$ whose coefficients tend to zero, so converging in the closed unit polydisk. It then makes sense to evaluate them at all ($n$-tuples of) roots of unity. She then proves:

\begin{theorem}
\label{thm:Ana}
\cite{Nei02}
If $f(x_1,\ldots,x_n) \in R[[x_1,\ldots,x_n]]$ converges in the closed unit polydisk, there exists $c>0$ such that, for any roots of unity $\zeta_1,\ldots,\zeta_n \in K$, either $F(\zeta_1,\ldots,\zeta_n) =0$ or $|F(\zeta_1,\ldots,\zeta_n)| \ge c$.
\end{theorem}

Prior to all this work, P. Monsky \cite{03681880} considered the special case of $p$-th power roots of unity and, subsequently, V. Serban \cite{07081553} has taken up this case.
Monsky and Serban instead consider arbitrary power series in $R[[x_1,\ldots,x_n]]$. These power series converge in the unit open polydisk and can be evaluated on $n$-tuples of the form $(\zeta_1-1,\ldots,\zeta_n-1)$ where the $\zeta_i$ are roots of unity of $p$-power order. Specifically:

\begin{theorem}
\label{thm:p-th}
\cite{03681880},\cite{07081553}
If $f(x_1,\ldots,x_n) \in R[[x_1,\ldots,x_n]]$, there exists $c>0$ such that, for any roots of unity $\zeta_1,\ldots,\zeta_n \in K$ of $p$-th power order, either \linebreak
$F(\zeta_1-1,\ldots,\zeta_n-1) =0$ or $|F(\zeta_1-1,\ldots,\zeta_n-1)| \ge c$.
\end{theorem}

\section{Semiabelian varieties}
\label{sec:AV}

\begin{conjecture}
\label{main}
Let $A/K$ be a semiabelian variety and $X$ a closed subvariety. 
There is  $c>0$ such that for any torsion point $P \in A(K)$ either $P \in X$ or $d(P,X) \ge c$.

\end{conjecture}

We stated this conjecture and proved it in the case $A$ is a torus and $K = \C_p$ \cite{zbMATH00972701} as a consequences of our result mentioned in Section \ref{sec:roots}.

After some partial results by A. Buium \cite{zbMATH00946959} and the author \cite{zbMATH00991915}, the conjecture was proved by T. Scanlon \cite{zbMATH01385059} when $K=\bar{\Q}_p$ and another proof given by C. Corpet \cite{arXiv:1309.7237}, also for $K=\bar{\Q}_p$. There is no compelling reason to restrict the conjecture to this case. We note that the equicharacteristic zero should be accessible, using the methods of \cite{zbMATH00537350}, for example. For the case of equicharacteristic $p>0$, R. de Jong, N. Looper and F. Shokrieh (in preparation) have shown that the statement of Conjecture \ref{main} holds if, in addition, $A$ and $X$ are defined over a global function field and $P$ has order prime to $p$. New techniques are needed in the case of $p$-th power torsion.

Let $A/K$ be a totally degenerate abelian variety in the sense of Mumford \cite{zbMATH03380769}. Then $A(K)$ is analytically isomorphc to $\G_m(K)^g/\Lambda$, where $\Lambda$ is a lattice in $\G_m(K)^g$. In particular, the torsion subgroup of $A(K)$ contains the image of the $g$-tuples of roots of unity and, using the results mentioned in Section \ref{sec:roots}, Conjecture \ref{main} can proved for this subgroup of the torsion subgroup. It would be very interesting to upgrade this argument to a proof of the full conjecture in this case, which would be an indication that no restriction on $K$ is necessary.

A version of Conjecture \ref{main} for Drinfeld modules was considered by D. Ghioca \cite{zbMATH05161936}.

\section{Siegel moduli spaces}

The first person to consider an analogue of Conjecture \ref{main} for Siegel moduli spaces (with CM points replacing torsion points) was P. Habegger \cite{zbMATH06340374} who proved a result for products of modular curves. 

For a fixed CM abelian variety $A/K$ with good ordinary reduction, the set of abelian varieties over $K$ with the same reduction as $A$ can described by the Serre-Tate parameters $\G_m(R)^{g(g+1)/2}$ and the set of special points in this set is formed by the points whose coordinates are the roots of unity of $p$-power order. Using the work of Monsky and Serban, this leads to a result about special points within this neigbourhood of $A$.

C. Qiu went much further and obtained a general result, using the theory of perfectoid spaces (which had been used earlier in the context of endomorphisms by J. Xie \cite{zbMATH06976300}) as follows.

\begin{theorem}
\label{thm:qiu}
\cite{zbMATH07683976}
Let $A/\bar{\Q}_p$ be a product of Siegel moduli spaces with arbitrary level structure and $X$ a closed subvariety.
There is  $c>0$ such that, for any ordinary CM point $P \in A(K)$, either $P \in X$ or $d(P,X) \ge c$.
\end{theorem}

Qiu also speculated about the possibility of a similar statement for more general Shimura varieties of Hodge type.

The restriction to special points of ordinary reduction is necessary. As shown by Habegger \cite{zbMATH06340374}, already in the case of elliptic curves, the set of CM points with supersingular reduction is not discrete.

\section{Algebraic Dynamics}

Buium \cite{zbMATH00946959} considered the periodic points of a lift of Frobenius on an algebraic abelian variety in connection with his results on Conjecture \ref{main} mentioned in Section \ref{sec:AV}. 

J. Xie \cite{zbMATH06976300} considered the backwards orbit of a point under a lift of Frobenius in $\PP^n/\C_p$ and introduced the idea of using the techniques of Scholze's perfectoid spaces. W. Peng \cite{zbMATH07679874} extended this to maps of $\A^n/\C_p$ which are lifts of the $p$-th power of another map. This means the following: If $F=(F_1,\ldots,F_n), \in R[x_1,\ldots,x_n]^n$, then $F$ is a lift of the $p$-th power of $G=(G_1,\ldots,G_n) \in R[x_1,\ldots,x_n]^n$
if $F_i = G_i^p$ in $R/M_{\e}$ for some $\e < 1$ (notation as in Section \ref{sec:prelim}). Peng's result is then:

\begin{theorem}
\label{thm:peng}
\cite{zbMATH07679874}
Let $F$ be a lift of a $p$-th power as above, viewed as a map $F:\A^n(\C_p) \to \A^n(\C_p)$ and $X$ a closed subvariety of $\A^n$. Then there exists  $c>0$ such that, for any periodic point $P \in \A^n(\C_p)$ for $F$, either $P \in X$ or $d(P,X) \ge c$.
\end{theorem}

\section{Loose ends}

The Archimedean case is quite different because the set of roots of unity in $\C$ is not discrete, in fact is dense in the unit circle. 

The question of how small a sum of a fixed number of complex roots of unity can be in terms of the maximum order of these roots of unity has been discussed in the literature, first by G. Myerson \cite{zbMATH03970845} and, later, at a lively discussion on a question of T. Tao on MathOverflow \cite{46068}. If $f(k,N)$ is the minimum nonzero complex absolute value of a sum of $k$ roots of unity of order dividing $N$, then $f(k,N) \ll N^{-\lambda_k}$, for some $\lambda_k > 0$ and one expects that this is best possible (for an appropriate choice of $\lambda_k$), although the best known lower bound on $f(k,N)$ is exponentially small in $N$. 

If one allows arbitrary coefficients on a linear form in complex roots of unity then one cannot expect much. For instance, if $\alpha$ is a Liouville number, then the absolute value of $e^{2\pi i\alpha}\zeta - 1$ can be made exponentially small in $N$ for a suitable choice of $N$-th root of unity $\zeta$, for infinitely many $N$. Maybe restricting the coefficients to algebraic numbers will be similar to the situation of sums of roots of unity although this seems to be completely unexplored, as are analogous questions about special points as above.

If one considers algebraic coefficients on a linear form in roots of unity $\sum a_i \zeta_i$, one can ask global questions about the behavior of $\max_{\zeta_i} |\sum a_i \zeta_i|_v$ as $v$ varies over the places of a ground field containing the $a_i$. In the function field case, one has a very nice global theory from classical algebraic geometry that even allows the $\zeta_i$ to vary among all constants. This was a main motivation for the work in \cite{zbMATH00972701} and, at the end of the paper, we speculated on how such a global theory might look in the number field case. There is, unfortunately, not much to report here except that a conjecture of A. Smirnov \cite {zbMATH00146392} (for the simplest case of a form in two variables) mentioned there was disproved by D. Masser \cite{zbMATH01766344} in the meanwhile.


\bibliographystyle{amsplain}
\bibliography{tate}

\end{document}